\documentclass[10pt]{article}
\usepackage{epsfig,amssymb,amsmath}
\usepackage{palatino}
\usepackage{verbatim}

\def\PP{{\cal P}}

\def\wt{\widetilde}

\def\la{\lambda}
\def\al{\alpha}

\newtheorem{lemma}{Lemma}[section]
\newtheorem{proposition}[lemma]{Proposition}
\newtheorem{corollary}[lemma]{Corollary}
\newtheorem{theorem}[lemma]{Theorem}
\newtheorem{remark}[lemma]{Remark}

\makeatletter
\@addtoreset{equation}{section}
\makeatother

\def\Frac#1#2{\mbox{\large${\textstyle \frac{#1}{#2}}$}}
\def\proof{\medskip\noindent{\bf Proof.} }
\def\qed{\hfill $\Box$}

\newcommand {\diag} {{\rm diag\,}}

\begin{document}

\title{Markov $L_2$ inequality with the Gegenbauer weight}

\author{Dragomir Aleksov, Geno Nikolov}

\date{}
\maketitle


\begin{abstract}

For the Gegenbauer weight function
$w_{\lambda}(t)=(1-t^2)^{\lambda-1/2}$, $\lambda>-1/2$, we denote by
$\Vert\cdot\Vert_{w_{\lambda}}$ the associated $L_2$-norm,
$$
\Vert
f\Vert_{w_{\lambda}}:=\Big(\int_{-1}^{1}w_{\lambda}(t)f^2(t)\,dt\Big)^{1/2}.
$$
We study the Markov inequality
$$
\Vert p^{\prime}\Vert_{w_{\lambda}}\leq c_{n}(\lambda)\,\Vert
p\Vert_{w_{\lambda}},\qquad p\in \PP_n,
$$
where $\PP_n$ is the class of algebraic polynomials of degree not
exceeding $n$. Upper and lower bounds for the best Markov constant
$c_{n}(\lambda)$ are obtained, which are valid for all $n\in
\mathbb{N}$ and $\la>-\frac{1}{2}$.
\end{abstract}


\section{Introduction and statement of the results}
Throughout this paper $\PP_n$ stands for the class of algebraic
polynomials of degree not exceeding $n$.

For the Gegenbauer weight function $w_{\la}(t)=(1-t^2)^{\la-1/2}$,
$\la>-1/2$, we denote by $\Vert\cdot\Vert_{w_\la}$ the associated
$L_2$-norm,
$$
\Vert
f\Vert_{w_\la}:=\Big(\int_{-1}^{1}w_{\la}(t)f^2(t)\,dt\Big)^{1/2}.
$$

Here we study the Markov inequality in this norm for the first
derivative of polynomials from $\PP_n$, in particular, we are
interested in the best Markov constant
$$
c_{n}(\la)=\sup_{\mathop{}^{p\in\PP_n}_{p\ne
0}}\frac{\Vert{p^{\prime}\Vert}_{w_\la}} {\Vert p\Vert_{w_\la}}.
$$

Let us start with a brief account of the known results.

In the case $\la=\frac{1}{2}$ (the case of a constant weight
function), E. Schmidt proved that
$$
c_n(1/2)=\frac{(2n+3)^2}{4\pi}\,\Big(1-\frac{\pi^2-3}{3(2n+3)^2}
+\frac{16R}{(2n+3)^4}\Big)^{-1}\,,\quad -6<R<13\,.
$$

Nikolov \cite{n03} studied two other particular cases, $\la=0,\,1$,
and proved the following two-sided estimates for the corresponding
Markov constants:
\begin{equation}\label{e1}
\begin{array}{l}
  0.472135 n^2 \le c_n(0) \le 0.478849 (n+2)^2\,, \\[0.5ex]
  0.248549 n^2 \le c_n(1) \le 0.256861 (n+\Frac{5}{2})^2\,.
\end{array}
\end{equation}

In \cite{ans16} we obtained an upper bound for $c_n(\la)$, which is
valid for all $n$ and $\la$:
$$
   c_n(\la) \le \frac{(n+1)(n+2\la +1)}{2\sqrt{2\la+1}}\,.
$$

This result has been improved in the recent paper \cite{ns17}, where
the following theorem was proved:\medskip

\noindent \textbf{Theorem A} {\it For all $\la > -\frac{1}{2}$ and
$n\ge 3$, the best constant $c_n(\la)$ in the Markov inequality
$$
   \|p_n'\|_{w_\la} \le c_n(\la) \|p_n\|_{w_\la}\,,\qquad p_n \in \PP_n\,,
$$
admits the estimates
\begin{eqnarray}
   \frac{1}{4}\frac{n^2(n+\la)^2}{(\la+1)(\la+2)} \;
 < & [c_n(\la)]^2 &
 < \;\frac{n(n+2\la+2)^3}{(\la+2)(\la+3)},
       \qquad\quad \la \ge 2\,; \label{e2} \\
   \frac{(n+\la)^2(n+2\la')^2}{(2\la+1)(2\la+5)} \;
 < & [c_n(\la)]^2 &
 < \;\frac{(n+\la+\la''+2)^4}{2(2\la+1)\sqrt{2\la+5}},
       \qquad \la > -\Frac{1}{2}\,, \label{e3}
\end{eqnarray}
where $\la' = \min\,\{0,\la\}$, $\,\la'' = \max\,\{0,\la\}$.}
\medskip

It has been also proved in \cite{ns17} that
$$
   [c_n(\la)]^2 \asymp \frac{1}{\la^2}n(n+2\la)^3\,,
$$
which shows that the upper bound in \eqref{e2} has the right order
in both $n$ and $\la$. The lower bound in \eqref{e2} is inferior to
the one in \eqref{e3}, it appears in \eqref{e2} just to indicate
that, roughly, for a fixed $\la$ and large $n$ the sharp Markov
constant is identified within a factor not exceeding two. Although
the upper bound in \eqref{e3} is not of the right order with respect
to $\la$, for moderate $\la$ (say, $\la\leq 25$) it is superior to
the one in \eqref{e2}.

In the present paper we prove two-sided estimates for $c_{n}(\la)$,
valid for all $\la>-1/2$, which are of the same nature as (and
slightly sharper than) those in \eqref{e3}. The approaches for their
derivation however are different. In \cite{ns17}, the results are
obtained through estimation of appropriate matrix norms. Here, we
identify the reciprocal of the squared best Markov constant as the
smallest zero of a related orthogonal polynomial, then exploit the
associated three-term recurrence relation to evaluate its lower
degree coefficients and eventually derive estimates for its smallest
zero. Let us mention that a similar relation between the best
constant in the $L_2$ Markov inequality with the Laguerre weight
function and the smallest zero of an orthogonal polynomial is given
in \cite[p. 85]{pd02}, and in \cite{ns16} we applied a similar
approach to obtain bounds for the best Markov constant in the
Laguerre case.

Our main result is the following theorem:
\begin{theorem}\label{t1.1}
For all $n\geq 3$ and for every $\la>-\Frac{1}{2}$, the best
constant $c_n(\la)$ in the Markov inequality
\begin{equation}\label{e1.1}
\Vert p^{\prime}\Vert_{w_\la}\leq
c_{n}(\la)\,\Vert p\Vert_{w_\la},\qquad p\in \PP_n,
\end{equation}
admits the estimates
\begin{equation}\label{e1.2}
\frac{(n+1)\big(n+\la+\frac{1}{2}\big)^2(n+2\la)}{(2\la+1)(2\la+5)}
\leq c_n^2(\la)\leq
\Frac{\big(n+\Frac{5}{4}\la+\Frac{9}{8}\big)^4}{2(2\la+1)\sqrt{2\la+5}}\,.
\end{equation}
\end{theorem}

By setting $\la = 0,1$ in \eqref{e1.2}, we obtain an improvement of
the upper bounds in \eqref{e1}, and combination with the lower
bounds in \eqref{e1} yields rather tight estimates.
\begin{corollary}\label{c1.2}
For the Chebyshev weights $w_0(x) = \frac{1}{\sqrt{1-x^2}}$ and
$w_1(x) = \sqrt{1-x^2}$, we have
\begin{eqnarray*}
  0.472135\, n^2 \le c_n(0) \le 0.472871\, \Big(n+\frac{9}{8}\Big)^2\,, &&\\
  0.248549\, n^2 \le c_n(1) \le 0.250987\, \Big(n+\frac{19}{8}\Big)^2\,. &&
\end{eqnarray*}
\end{corollary}

For the proof of Theorem~\ref{t1.1} we obtain separately estimates
for $c_n(\la)$ in the cases of even and odd $n$ (Theorems~\ref{t4.2}
and \ref{t4.4}). These estimates are slightly sharper than the ones
in Theorem~\ref{t1.1}, in particular, they yield the following
asymptotic inequalities:

\begin{corollary}\label{c1.3}
For every $n\geq 3$, there holds
\begin{equation}\label{e1.3}
\frac{(n+2)(n-1)n^2}{4}\leq\lim_{\la\rightarrow
-\frac{1}{2}}(2\la+1)\,c_n^2(\la)\leq\frac{n^2(n+1)^{2}}{4}\,.
\end{equation}
\end{corollary}

The paper is organised as follows. In Sect. 2 we show that the
reciprocal of the squared best Markov constant, $1/[c_n(\la)]^2$, is
equal to the smallest zero of an orthogonal polynomial of degree
$m=\lfloor\frac{n+1}{2}\rfloor$ (different in  the cases $n=2m$ and
$n=2m-1$), and we derive the three-term recurrence relation
satisfied by these orthogonal polynomials. Based on the three-term
recurrence relations, in Sect. 3 we evaluate and estimate the lowest
degree coefficients of the $m$-th orthogonal polynomial. In Sect. 4
we prove estimates for $c_n(\la)$ in the cases of even and odd $n$
(Theorems~\ref{t4.2} and \ref{t4.4}), and derive as consequences
Theorem \ref{t1.1} and Corollary~\ref{c1.3}.
\section{\boldmath{$c_{n}^2(\la)$} and the
extreme zero of an orthogonal polynomial}

In a recent paper \cite{ans16} we showed that the extreme polynomial
in the Markov inequality \eqref{e1.1} is even or odd if $n$ is even
or odd. The following theorem summarizes some of the results
obtained in \cite{ans16}:
\begin{theorem}\label{t2.1}
The best constant $c_{n}(\la)$ in the Markov inequality \eqref{e1.1}
is given by
\begin{equation}\label{e2.1}
c_n(\la)=\begin{cases} 2\sqrt{\nu_m}\,, & \ n=2m, \vspace*{2mm}\\
2\sqrt{\wt{\nu}_m}\,, & \ n=2m-1\,,
\end{cases}
\end{equation}
where $\nu_m$ and $\wt{\nu}_m$ are the largest eigenvalues of the
$m\times m$ positive definite matrices
$\,\mathbf{C}_m^{\top}\mathbf{C}_m\,$ and
$\,\wt{\mathbf{C}}_m^{\top}\wt{\mathbf{C}}_m\,$, respectively, given
by
\begin{equation}\label{e2.2}
\mathbf{C}_m=\begin{pmatrix} \alpha_1\beta_1&\alpha_1\beta_2&\cdot\cdots&\alpha_1\beta_m\\
0&\alpha_2\beta_2&\cdots&\alpha_2\beta_m\\
\vdots&\vdots&\ddots&\vdots\\0&0&\cdots&\alpha_m\beta_m\end{pmatrix}\,,\qquad
\mathbf{\widetilde{C}}_m=
\begin{pmatrix}
\widetilde{\alpha}_1\widetilde{\beta}_1&\widetilde{\alpha}_1\widetilde{\beta}_2
&\cdot\cdots&\widetilde{\alpha}_1\widetilde{\beta}_m\\
0&\widetilde{\alpha}_2\widetilde{\beta}_2&\cdots&\widetilde{\alpha}_2\widetilde{\beta}_m\\
\vdots&\vdots&\ddots&\vdots\\0&0&\cdots&\widetilde{\alpha}_m\widetilde{\beta}_m\end{pmatrix}\,.
\end{equation}
Here,
\begin{eqnarray}
   \al_k := (2k-1 + \la) h_{2k-1},
& \beta_k := \Frac{1}{h_{2k}}\,; \label{e2.3}\\
   \wt\al_k := (2k-2 + \la) h_{2k-2},
& \wt\beta_k := \Frac{1}{h_{2k-1}}\,,\label{e2.4}
\end{eqnarray}
with
\begin{equation}\label{e2.5}
   h_i^2
:= h_{i,\lambda}^2 :=
\frac{\Gamma(i+2\lambda)}{(i+\lambda)\Gamma(i+1)}\,.
\end{equation}
\end{theorem}

Clearly, matrices $\mathbf{C}_m$ and $\wt{\mathbf{C}}_m$ can be
represented as
\begin{equation}\label{e2.7}
\begin{split}
&\mathbf{C}_m=\diag (\al_1,\ldots,\al_m)\mathbf{T}_m \diag
(\beta_1,\ldots,\beta_m)\,,\\
& \wt{\mathbf{C}}_m= \diag
(\wt{\al}_1,\ldots,\wt{\al}_m)\mathbf{T}_m\diag
(\wt{\beta}_1,\ldots,\wt{\beta}_m)\,,
\end{split}
\end{equation}
where $\mathbf{T}_m$ is an upper tri-diagonal $m\times m$ matrix
with non-zero entries equal to $1$,
$$
\mathbf{T}_m=\begin{pmatrix} 1 & 1 & \cdots & 1\\ 0 & 1 & \cdots & 1\\
\vdots & \vdots & \ddots & \vdots\\ 0 & 0 & \cdots & 1
\end{pmatrix}\,.
$$

Since
$\,\mathbf{C}_m^{\top}\mathbf{C}_m\sim\mathbf{C}_m\mathbf{C}_m^{\top}\,$
and $\,\wt{\mathbf{C}}_m^{\top}\wt{\mathbf{C}}_m\sim
\wt{\mathbf{C}}_m\wt{\mathbf{C}}_m^{\top}\,$, we conclude that
\begin{equation}\label{e2.8}
\begin{split}
&\nu_m\ \text{ is the largest eigenvalue of the matrix }\
\mathbf{A}_m:=\mathbf{C}_m\mathbf{C}_m^{\top}\,,\\
&\wt{\nu}_m\ \text{ is the largest eigenvalue of the matrix }\
\wt{\mathbf{A}}_m:=\wt{\mathbf{C}}_m\wt{\mathbf{C}}_m^{\top}\,.
\end{split}
\end{equation}

It turns out that it is advantageous to work with the inverse
matrices $\mathbf{B}_m:=\mathbf{A}_m^{-1}$ and
$\wt{\mathbf{B}}_m:=\wt{\mathbf{A}}_m^{-1}$, respectively, as
$\mathbf{B}_m$ and $\wt{\mathbf{B}}_m$ are tri-diagonal matrices.
Below we demonstrate this for $\mathbf{B}_m$.

The matrix $\mathbf{T}_m^{-1}$ is two-diagonal, namely,
\begin{equation}\label{e2.9}
\mathbf{T}_m^{-1}=\begin{pmatrix} 1 & -1 & 0 & \cdots & 0 \\
0 & 1 & -1 & \cdots & 0 \\ 0 & 0 & 1 & \cdots & 0 \\
\vdots & \vdots & \vdots & \ddots & \vdots \\
0 & 0 & 0 & \cdots & 1
\end{pmatrix}\,.
\end{equation}
For
$\mathbf{B}_m=\mathbf{A}_m^{-1}
=\big(\mathbf{C}_m^{\top}\big)^{-1}\mathbf{C}_m^{-1}
=\big(\mathbf{C}_m^{-1}\big)^{\top}\mathbf{C}_m^{-1}$, using
\eqref{e2.7}, we have
\[
\begin{split}
\mathbf{B}_m=&\big(\diag
(\beta_1^{-1},\ldots,\beta_m^{-1})\mathbf{T}_m^{-1}\diag
(\al_1^{-1},\ldots,\al_m^{-1})\big)^{\top} \diag
(\beta_1^{-1},\ldots,\beta_m^{-1})\mathbf{T}_m^{-1}\diag
(\al_1^{-1},\ldots,\al_m^{-1})\\
=&\diag
(\al_1^{-1},\ldots,\al_m^{-1})\big(\mathbf{T}_m^{-1}\big)^{\top}
\diag (\beta_1^{-2},\ldots,\beta_m^{-2})\mathbf{T}_m^{-1} \diag
(\al_1^{-1},\ldots,\al_m^{-1})\,.
\end{split}
\]
Making use of \eqref{e2.9}, we perform the multiplications to
conclude that, indeed, $\mathbf{B}_m$ is tri-diagonal.  We formulate
the result below:
\begin{proposition}\label{p2.2}
The matrix $\mathbf{A}_m^{-1}=:\mathbf{B}_m=(b_{i,j})_{m\times m}$
is symmetric and tri-diagonal, with elements
\begin{eqnarray}
&&b_{1,1}=  \Frac{1}{\al_1^2\beta_1^2}\,, \label{e2.10}\\
&&b_{k,k}=\Frac{1}{\al_{k}^2}\Big(\Frac{1}{\beta_{k-1}^2}+\Frac{1}{\beta_{k}^2}\Big)\,,
\quad k=2,\ldots,m\,,\label{e2.11}\\
&&b_{k,k+1}=-\Frac{1}{\al_k\al_{k+1}\beta_{k}^2}\,,\quad
k=1,\ldots,m-1\,.\label{e2.12}
\end{eqnarray}

The same conclusion applies to the matrix
$\,\wt{\mathbf{A}}_m^{-1}=:\wt{\mathbf{B}}_m=(\wt{b}_{i,j})_{m\times
m}$, with the $\,b$'s, $\,\al$'s and $\,\beta$'s replaced by the
$\,\wt{b}$'s, $\,\wt{\al}$'s and $\,\wt{\beta}$'s.
\end{proposition}

Thus, $\mathbf{B}_m$ and $\wt{\mathbf{B}}_m$ are Jacobi matrices,
which are positive definite as inverse of the positive definite
matrices $\mathbf{A}_m$ and $\wt{\mathbf{A}}_m$. The characteristic
polynomials of $\mathbf{B}_m$ and $\wt{\mathbf{B}}_m$,
$$
P_m(\mu)=\det(\mu \mathbf{E}_m-\mathbf{B}_m)\,,\qquad
\wt{P}_m(\mu)=\det(\mu \mathbf{E}_m-\wt{\mathbf{B}}_m),
$$
are determined by three-term recurrence relations, and, by Favard's
theorem, $\{P_m\}$ and $\{\wt{P}_m\}$ constitute two sequences of
orthogonal polynomials with respect to measures supported on the
positive axis. Let $\mu_1<\mu_2<\cdots<\mu_m$ and
$\wt{\mu}_1<\wt{\mu}_2<\cdots<\wt{\mu}_m$ be the zeros of $P_m$ and
$\wt{P}_m$, respectively, i.e., the eigenvalues of $\mathbf{B}_m$
and $\wt{\mathbf{B}}_m$. Since the latter are reciprocal to the
eigenvalues of $\mathbf{A}_m$ and $\wt{\mathbf{A}}_m$, in
particular, $\nu_m=\mu_1^{-1}$ and $\wt{\nu}_m=\wt{\mu}_1^{-1}$,
Theorem~\ref{t2.1}, \eqref{e2.8} and Proposition~\ref{p2.2} yield
the following
\begin{theorem}\label{t2.3}
The best constant $c_{n}(\la)$ in the Markov inequality \eqref{e1.1}
is given by
\begin{equation}\label{e2.13}
c_n(\la)=\begin{cases} \Frac{2}{\sqrt{\mu_1}}\,, & \ n=2m, \vspace*{2mm}\\
\Frac{2}{\sqrt{\wt{\mu}_1}}\,, & \ n=2m-1\,,
\end{cases}
\end{equation}
where $\mu_1$ and $\wt{\mu}_1$ are the smallest zeros of monic
polynomials $P_m$ and $\wt{P}_m$, orthogonal with respect to a
measure supported on $\mathbb{R}_{+}$. The polynomials $\{P_k\}$ are
defined by the three-term recurrence relation
\begin{equation}\label{e2.14}
\begin{split}
&P_k(\mu)=\Big[\mu-\frac{1}{\al_k^2}\Big(\frac{1}
{\beta_{k-1}^2}+\frac{1}{\beta_k^2}\Big)\Big]P_{k-1}(\mu)
-\frac{1}{\al_{k-1}^2\al_k^2\beta_{k-1}^4}P_{k-2}(\mu)\,,\quad k\geq
2\,,\\
& P_{0}(\mu)=1\,,\quad P_{1}(\mu)=\mu-\frac{1}{\al_1^2\beta_1^2}\,.
\end{split}
\end{equation}
The polynomials $\{\wt{P}_k\}$ satisfy the same recurrence relation,
with the $\al$'s and $\beta$'s replaced by the $\wt{\al}$'s and
$\wt{\beta}$'s.
\end{theorem}

We renormalise polynomials $\{P_k\}_{0}^{m}$ and
$\{\wt{P}_k\}_{0}^{m}$ by setting $Q_0=P_0=\wt{Q}_0=\wt{P}_0=1$ and
$$
Q_k=d_k\,P_k\,,\quad \wt{Q}_k=\wt{d}_k\,\wt{P}_k\,,\quad
k=1,\ldots,m
$$
so that
\begin{equation}\label{e2.15}
Q_k(0)=\wt{Q}_k(0)=1,\qquad k=0,\ldots,m
\end{equation}
(note that this is possible because all the zeros of $P_k$ and
$\wt{P}_k$ are positive).

For $k=1,\ldots,m$, we have $P_k(\mu)=\det(\mu
\mathbf{E}_k-\mathbf{B}_k)$ and $\wt{P}_k(\mu)=\det(\mu
\mathbf{E}_k-\wt{\mathbf{B}}_k)$, therefore,
\begin{eqnarray*}
&&P_k(0)=\det(-\mathbf{B}_k)=(-1)^k\det(\mathbf{B}_k)
=(-1)^k\det(\mathbf{A}_k^{-1})
=(-1)^k\det(\mathbf{A}_k)^{-1}\,,\\
&&\wt{P}_k(0)=\det(-\wt{\mathbf{B}}_k)=(-1)^k\det(\wt{\mathbf{B}}_k)
=(-1)^k\det(\wt{\mathbf{A}}_k^{-1})
=(-1)^k\det(\wt{\mathbf{A}}_k)^{-1}\,.
\end{eqnarray*}
Since $\mathbf{A}_k=\mathbf{C}_k\mathbf{C}_k^{\top}$ and
$\wt{\mathbf{A}}_k=\wt{\mathbf{C}}_k\wt{\mathbf{C}}_k^{\top}$\,, we
make use of \eqref{e2.7} (with $m$ replaced by $k$) to obtain
$$
\det(\mathbf{A}_k)=\det(\mathbf{C}_k)^2
=\prod_{i=1}^{k}\al_i^2\beta_i^2\,,\quad
\det(\wt{\mathbf{A}}_k)=\det(\wt{\mathbf{C}}_k)^2
=\prod_{i=1}^{k}\wt{\al}_i^2\wt{\beta}_i^2\,.
$$
Consequently,
$$
P_k(0)=\frac{(-1)^k}{\prod_{i=1}^{k}\al_i^2\beta_i^2}\,,\quad
\wt{P}_k(0)=\frac{(-1)^k}{\prod_{i=1}^{k}\wt{\al}_i^2\wt{\beta}_i^2}\
\Rightarrow\ d_k=(-1)^k\,\prod_{i=1}^{k}\al_i^2\beta_i^2\,\quad
\wt{d}_k=(-1)^k \prod_{i=1}^{k}\wt{\al}_i^2\wt{\beta}_i^2\,.
$$
Thus, the renormalised to satisfy \eqref{e2.15} polynomials
$\{Q_{k}\}$ and $\{\wt{Q}_k\}$ are given by
\begin{equation}\label{e2.16}
Q_k(\mu)=(-1)^k\,\Big(\prod_{i=1}^{k}\al_i^2\beta_i^2\Big)\,P_k(\mu)\,,\qquad
\wt{Q}_k(\mu)=(-1)^k
\Big(\prod_{i=1}^{k}\wt{\al}_i^2\wt{\beta}_i^2\Big)\,\wt{P}_k(\mu)\,,
\quad k=1,\ldots,m.
\end{equation}

From \eqref{e2.14} it is easy to deduce the recurrence relations
satisfied by $\{Q_{k}\}$ and $\{\wt{Q}_k\}$.

\begin{proposition}\label{p2.4}
The polynomials $\{Q_{k}\}$ in \eqref{e2.16} satisfy the recurrence
relation
\begin{equation}\label{e2.17}
\begin{split}
&Q_k(\mu)-Q_{k-1}(\mu)=\Frac{\beta_k^2}{\beta_{k-1}^2}\,
\big[Q_{k-1}(\mu)-Q_{k-2}(\mu)\big]-\al_k^2\beta_k^2\,\mu\,Q_{k-1}(\mu)\,,
\qquad k\geq 2\,,\\
&Q_{0}(\mu)=1\,,\quad Q_1(\mu)=1-\al_1^2\beta_1^2\,\mu\,.
\end{split}
\end{equation}

The polynomials $\{\wt{Q}_{k}\}$ in \eqref{e2.16} satisfy the same
recurrence relation, with the $\al$'s and $\beta$'s replaced by the
$\wt{\al}$'s and $\wt{\beta}$'s.
\end{proposition}

\section{The lowest degree coefficients of \boldmath{$Q_m$}
and \boldmath{$\wt{Q}_m$}}

In view of \eqref{e2.15}, we may write polynomials $Q_k$ and
$\wt{Q}_k$, $k\geq 1$, in the form
\begin{equation} \label{e3.1}
\begin{split}
& Q_k(\mu)=1-A_{1,k}\,\mu+A_{2,k}\,\mu^2-\cdots+(-1)^{k}A_{k,k}\,\mu^k\,,\\
& \wt{Q}_k(\mu)=1-\wt{A}_{1,k}\,\mu+\wt{A}_{2,k}\,\mu^2+\cdots
+(-1)^{k}\wt{A}_{k,k}\,\mu^k \,.
\end{split}
\end{equation}

Our goal now is to find expressions for $A_{i,m},\;\wt{A}_{i,m}$,
$i=1,2$. First of all, we make use of \eqref{e2.3}--\eqref{e2.4} to
find the explicit form of the coefficients occurring in recurrence
formulae for $Q_k$ and $\wt{Q}_k$. We have
\begin{eqnarray}
&& \Frac{\beta_k^2}{\beta_{k-1}^2}=\Frac{k(2k-1)(2k+\la)}
{(k-1+\la)(2k-2+\la)(2k-1+2\la)}\,,\qquad
\al_k^2\beta_k^2=\Frac{2k(2k-1+\la)(2k+\la)}{2k-1+2\la}\,,
\label{e3.2}\vspace*{2mm}\\
&&
\Frac{\wt{\beta}_k^2}{\wt{\beta}_{k-1}^2}=\Frac{(k-1)(2k-1)(2k-1+\la)}
{(k-1+\la)(2k-3+\la)(2k-3+2\la)}\,,\qquad
\wt{\al}_k^2\wt{\beta}_k^2=\Frac{(2k-1)(2k-2+\la)(2k-1+\la)}{2(k-1+\la)}\,.
\label{e3.3}
\end{eqnarray}
By substituting these quantities in the recurrence formulae  in
Proposition \ref{p2.4} and replacing $k$ by $m$, we obtain
\begin{equation}\label{e3.4}
\begin{split}
Q_m(\mu)-Q_{m-1}(\mu)=&\frac{m(2m-1)(2m+\la)}{(m-1+\la)(2m-2+\la)(2m-1+2\la)}\,
\big[Q_{m-1}(\mu)-Q_{m-2}(\mu)\big]\\
&-\frac{2m(2m-1+\la)(2m+\la)}{2m-1+2\la}\,\mu\,Q_{m-1}(\mu)\,,
\end{split}
\end{equation}
\begin{equation}\label{e3.5}
\begin{split}
\wt{Q}_m(\mu)-\wt{Q}_{m-1}(\mu)=&
\frac{(m-1)(2m-1)(2m-1+\la)}{(m-1+\la)(2m-3+\la)(2m-3+2\la)}\,
\big[\wt{Q}_{m-1}(\mu)-\wt{Q}_{m-2}(\mu)\big]\\
&-\frac{(2m-1)(2m-2+\la)(2m-1+\la)}{2(m-1+\la)}\,\mu\,\wt{Q}_{m-1}(\mu)\,.
\end{split}
\end{equation}

\begin{lemma}\label{l3.1}
For every $m\in \mathbb{N}_0$ there holds
$$
(i)\ \ A_{1,m}=\frac{m(m+1)(m+\la)(m+\la+1)}{2\la+1}\,;\qquad (ii) \
\
\wt{A}_{1,m}=\frac{m(m+\la)\big(m^2+\la\,m-\Frac{1}{2}\big)}{2\la+1}\,.
$$
\end{lemma}

\proof  (i) The formula is true for $m=0$, since $Q_0(\mu)=1$, and
hence $A_{1,0}=0$. Clearly, (i) holds for $m=1$, too, since, by
\eqref{e2.17} and \eqref{e3.2},
$$
A_{1,1}=\al_1^2\beta_1^2=\frac{2(\la+1)(\la+2)}{2\la+1}\,.
$$
We set $D_{1,k}:=A_{1,k}-A_{1,k-1}$, $k\in \mathbb{N}$, then claim
(i) is equivalent to
\begin{equation}\label{e3.6}
D_{1,m}=\frac{2m(m+\la)(2m+\la)}{2\la+1}\,,\quad m\in \mathbb{N}\,,
\end{equation}
and it is true for $m=1$, since $D_{1,1}=A_{1,1}$. We shall prove
\eqref{e3.6} by induction with respect to $m$. To this end, we
differentiate \eqref{e3.4} in $\mu$ and then set $\mu=0$, making use
of \eqref{e3.1}, to obtain the recurrence formula
$$
D_{1,m}=\frac{m(2m-1)(2m+\la)}{(m-1+\la)(2m-2+\la)(2m-1+2\la)}\,D_{1,m-1}
+\frac{2m(2m-1+\la)(2m+\la)}{2m-1+2\la}\,.
$$
Assuming that \eqref{e3.6} is true for $m-1$, $m\geq 2$, we
substitute the expression for $D_{m-1}$ in the above formula to
verify that \eqref{e3.6} holds for $m$:
\[
\begin{split}
D_{1,m}&=\frac{m(2m-1)(2m+\la)2(m-1)(m-1+\la)(2m-2+\la)}
{(m-1+\la)(2m-2+\la)(2m-1+2\la)(2\la+1)}
+\frac{2m(2m-1+\la)(2m+\la)}{2m-1+2\la}\vspace*{2mm}\\
&=\frac{2m(m+\la)(2m+\la)}{2\la+1}\,.
\end{split}
\]

(ii) Clearly, (ii) holds for $m=0$, since $\wt{A}_{1,0}=0$, and it
is also true for $m=1$, since, by Proposition~\ref{p2.4} and
\eqref{e3.3},
$$
\wt{A}_{1,1}=\wt{\al}_1^2\wt{\beta}_1^2=\frac{\la+1}{2}\,.
$$
Similarly to the proof of (i), we set
$\wt{D}_{1,k}=\wt{A}_{1,k}-\wt{A}_{1,k-1}$, $k\in \mathbb{N}$, then
(ii) is equivalent to
\begin{equation}\label{e3.7}
\wt{D}_{1,m}=\frac{(2m-1)(2m-1+\la)(2m-1+2\la)}{2(2\la+1)}\,,
\end{equation}
and the latter is true for $m=1$, since $\wt{D}_{1,1}=\wt{A}_{1,1}$.
Similarly to the proof of (i), we obtain a recurrence relation by
differentiating \eqref{e3.5} and then substituting $\mu=0$:
$$
\wt{D}_{1,m}=\frac{(m-1)(2m-1)(2m-1+\la)}
{(m-1+\la)(2m-3+\la)(2m-3+2\la)}\,\wt{D}_{1,m-1}
+\frac{(2m-1)(2m-2+\la)(2m-1+\la)}{2(m-1+\la)}\,.
$$

We observe that the right-hand side of \eqref{e3.7} is obtained from
the right-hand side of \eqref{e3.6} by the change $m\mapsto m-1/2$,
and the same change transforms the recurrence relation for $D_m$
into the recurrence relation for $\wt{D}_m$. Therefore, \eqref{e3.7}
is a consequence of \eqref{e3.6}. \qed

\begin{remark}\label{r3.2}
The coefficients $A_{1,m}$ and $\wt{A}_m$ are in fact the traces of
matrices $\mathbf{A}_m$ and $\wt{\mathbf{A}}_m$, respectively, and
they were evaluated in \cite[Lemma~2.3]{ans16}. We incorporate an
alternative proof first, for the sake of completeness and, second,
because the same approach is applied below for the evaluation of
coefficients $A_{2,m}$ and $\wt{A}_{2,m}$.
\end{remark}

Next, we proceed with the evaluation of the coefficients $A_{2,m}$
and $\wt{A}_{2,m}$. Let us set
\begin{eqnarray}
&&D_{2,1}=0\,,\quad D_{2,m}:=A_{2,m}-A_{2,m-1}\,,\quad m\geq 2\,,
\label{e3.8}\\
&&\wt{D}_{2,1}=0\,,\quad
\wt{D}_{2,m}:=\wt{A}_{2,m}-\wt{A}_{2,m-1}\,,\quad m\geq 2\,.
\label{e3.9}
\end{eqnarray}

\begin{lemma}\label{l3.2}
(i) The sequence $\{D_{2,m}\}$ defined by \eqref{e3.8} satisfies the
recurrence relation
\begin{equation}\label{e3.10}
\begin{split}
D_{2,m}=&\frac{m(2m-1)(2m+\la)}{(m-1+\la)(2m-2+\la)(2m-1+2\la)}\,D_{2,m-1}\\
&+\frac{2(m-1)m^2(m-1+\la)(m+\la)(2m-1+\la)(2m+\la)}{(2\la+1)(2m-1+2\la)}\,.
\end{split}
\end{equation}

The solution of \eqref{e3.10} with the initial condition $D_{2,1}=0$
is given by
\begin{equation}\label{e3.11}
D_{2,m}=\frac{2(m-1)m(m+\la)(m+\la+1)(2m+\la)\big[m^2+\la
m-\frac{2}{2\la+3}\big]} {(2\la+1)(2\la+5)}\,.
\end{equation}

(ii) The sequence $\{\wt{D}_{2,m}\}$ defined by \eqref{e3.9}
satisfies the recurrence relation
\begin{equation}\label{e3.12}
\begin{split}
\wt{D}_{2,m}=&\frac{(m-1)(2m-1)(2m-1+\la)}
{(m-1+\la)(2m-3+\la)(2m-3+2\la)}\,\wt{D}_{2,m-1}\\
&+\frac{(m-1)(2m-1)(2m-2+\la)(2m-1+\la)
\big[m^2+(\la-2)m-\la+\frac{1}{2}\big]}{2(2\la+1)}\,.
\end{split}
\end{equation}

The solution of \eqref{e3.12} with the initial condition
$\wt{D}_{2,1}=0$ is given by
\begin{equation}\label{e3.13}
\wt{D}_{2,m}=\frac{(m\!-\!1)(2m\!-\!1)(m\!+\!\la)
(2m\!-\!1\!+\!\la)(2m\!-\!1\!+\!2\la)
\big[m^2\!+\!(\la\!-\!1)m\!-\!\frac{2\la\!+1}{2}\!-\!
\frac{2}{2\la\!+\!3}\big]} {2(2\la+1)(2\la+5)}\,.
\end{equation}
\end{lemma}

\proof The recurrence formula \eqref{e3.10} is deduced by two-fold
differentiation of \eqref{e3.4} with respect to $\mu$, then setting
$\mu=0$ and using Lemma~\ref{l3.1}(i) to replace $A_{1,m-1}$ in the
resulting identity. The recurrence formula \eqref{e3.12} is obtained
in the same manner: we differentiate \eqref{e3.5} twice, then set
$\mu=0$ and apply Lemma~\ref{l3.1}(ii) to replace $\wt{A}_{1,m-1}$
in the resulting identity.

Now it is a straightforward (tough rather tedious) task to verify
that the sequences $\{D_{2,m}\}$ and $\{\wt{D}_{2,m}\}$ defined by
\eqref{e3.11} and \eqref{e3.13} are the solutions of the recurrence
relations \eqref{e3.10} and \eqref{e3.12}, respectively, with the
initial conditions $D_{2,1}=0$, $\wt{D}_{2,1}=0$. \qed

\begin{lemma}\label{l3.3}
The coefficients $A_{2,m}$ and $\wt{A}_{2,m}$ are given by
\begin{equation}\label{e3.14}
A_{2,m}=\frac{(m-1)m(m+1)(m+\la)(m+\la+1)(m+\la+2)
\big[m^2+(\la+1)m+\Frac{4\la^2+2\la-14}{3(2\la+3)}\big]}
{2(2\la+1)(2\la+5)}
\end{equation}
and
\begin{eqnarray}
&&\wt{A}_{2,m}=\frac{(m-1)m(m+\la)(m+\la+1)r_{\la}(m) }
{24(2\la+1)(2\la+3)(2\la+5)}\,,\vspace*{2mm} \label{e3.15}\\
&& \begin{array}{l}r_{\la}(m):=12(2\la+3)m^4+24\la(2\la+3)m^3
+4(6\la^3+7\la^2-19\la-32)m^2\\
\qquad\qquad\quad -4\la(2\la^2+19\la+32)m-8\la^3-20\la^2+14\la+71\,.
\end{array}
\label{e3.16}
\end{eqnarray}
\end{lemma}

\proof We have
$$
A_{2,m}=\sum_{k=2}^m D_{2,k}\,,\qquad \wt{A}_{2,m}=\sum_{k=2}^m
\wt{D}_{2,k}\,,
$$
hence, knowing formulae \eqref{e3.14} and
\eqref{e3.15}--\eqref{e3.16}, one may think of proving them by
induction with respect to $m$, especially having in mind that the
induction base is obvious. However, performing the induction step by
hand, though possible, is a hard work, this is why  we highly
recommend for that purpose the usage of a computer algebra program,
for instance, Wolfram's Mathematica does perfectly that job.

A reasonable question here is: how do we guess formulae
\eqref{e3.14} and \eqref{e3.15}--\eqref{e3.16}? Our approach makes
use of the observation that $A_{2,m}$ and $\wt{A}_{2,m}$ are
polynomials in $m$. We evaluate these coefficients for several
consecutive values of $m$ (nine values suffice!) and then construct
the associated interpolating polynomials to deduce the expressions
for $A_{2,m}$ and $\wt{A}_{2,m}$. Needles to say, we have used a
computer algebra program for this purpose. \qed

Next, we obtain two-sided estimates for the coefficients $A_{2,m}$
and $\wt{A}_{2,m}$.

\begin{lemma}\label{l3.4}
For all $m\in \mathbb{N}$, $m\geq 2$, and for every
$\la>-\Frac{1}{2}$, the coefficient $A_{2,m}$ admits the estimates
$$
\frac{(m\!-\!1)m^2(m\!+\!1)(m\!+\!\la)^2(m\!+\!\la\!+\!1)^2}
{2(2\la+1)(2\la+5)} \leq A_{2,m}\leq
\frac{(m\!-\!1)m(m\!+\!1)^2(m\!+\!\la)^2(m\!+\!\la\!+\!1)(m\!+\!\la\!+\!2)}
{2(2\la+1)(2\la+5)}\,.
$$
\end{lemma}

\proof We use formula \eqref{e3.14}. For the lower estimate, we need
to show that
$$
(m+\la+2)\Big[m^2+(\la+1)m
+\frac{4\la^2+2\la-14}{3(2\la+3)}\Big]\geq m(m+\la)(m+\la+1)\,.
$$
The difference of the left-hand and the right-hand sides is equal to
$$
g_0(m):=2m^2+\frac{4(4\la^2+8\la+1)}{3(2\la+3)}\,m
+\frac{2(\la+2)(2\la^2+\la-7)}{3(2\la+3)}\,.
$$
It is easy to see that $g_0^{\prime}(m)>0$ for $m\geq 2$ and
$\la>-1/2$, therefore $g_0$ is monotone increasing, and
$$
g_0(m)\geq g_0(2)=\frac{2(2\la^3+21\la^2+51\la+26)}{3(2\la+3)}>0\,.
$$

For the upper estimate, we need to prove the inequality
$$
m^2+(\la+1)m +\frac{4\la^2+2\la-14}{3(2\la+3)}\le (m+1)(m+\la)\,.
$$
The latter is equivalent to the inequality
$$
\frac{4\la^2+2\la-14}{3(2\la+3)}<\lambda\,,
$$
which is readily verified to be true for $ \lambda>-1/2$. \qed

\begin{lemma}\label{l3.5}
For all $m\in \mathbb{N}$, $m\geq 2$, the coefficient $\wt{A}_{2,m}$
admits the lower estimates
\begin{eqnarray*}
&&(i)\quad \wt{A}_{2,m}\geq \frac{(m-1)m(m+\la)(m+\la+1)\big(m^2+\la
m-\frac{1}{2}\big)\big(m^2+\la m-\frac{\la}{3}-\frac{7}{2}\big)}
{2(2\la+1)(2\la+5)}\,,\quad -\frac{1}{2}<\la\leq 0\,,\vspace*{2mm}\\
&&(ii)\quad \wt{A}_{2,m}\geq
\frac{(m-1)m(m+\la)(m+\la+1)\big(m^2+\la
m-\frac{1}{2}\big)\big(m^2+\la m-\frac{\la}{2}-\frac{7}{2}\big)}
{2(2\la+1)(2\la+5)}\,,\quad \la\geq 0\,.
\end{eqnarray*}

For all $m\in \mathbb{N}$, $m\geq 2$, and for every
$\la>-\Frac{1}{2}$, the coefficient $\wt{A}_{2,m}$ admits the upper
estimate
$$
\wt{A}_{2,m}\leq\frac{(m-1)m^2(m+\la)^2(m+\la+1)\big(m^2+\la
m-\frac{1}{2}\big)}{2(2\la+1)(2\la+5)}\,.
$$
\end{lemma}

\proof The polynomial $r_{\la}$ in \eqref{e3.16} satisfies
\begin{equation}\label{e3.17}
\begin{split}
r_{\la}(m)=&(2\la\!+\!3)\big(12m^4\!+\!24\la
m^3\!+\!(12\la^2\!-\!4\la\!-\!32)m^2
\!-\!(4\la^2\!+\!32\la\!+\!16)m\!-\!4\la^2\!-\!4\la\!+\!13\big)\\
&-16(m-2)(2m+1)\,,
\end{split}
\end{equation}
therefore
\[
\begin{split}
r_{\la}(m)&\leq (2\la\!+\!3)\big(12m^4\!+\!24\la
m^3\!+\!(12\la^2\!-\!4\la\!-\!32)m^2
\!-\!(4\la^2\!+\!32\la\!+\!16)m\!-\!4\la^2\!-\!4\la\!+\!13\big)\\
&=:(2\la+3)s_{\la}(m)\,.
\end{split}
\]

On the other hand,
\[
\begin{split}
s_{\la}(m)&=(12m^2+12\la m-4\la-26)\big(m^2+\la m-\frac{1}{2}\big)-
(16m+4\la^2+6\la)\\
&<12m(m+\la)\big(m^2+\la m-\frac{1}{2}\big)\,,
\end{split}
\]
hence
$$
r_{\la}(m)\leq 12(2\la+3)m(m+\la)\big(m^2+\la m-\frac{1}{2}\big)\,.
$$
The upper estimate for $\wt{A}_{2,m}$ now follows by putting this
upper bound for $r_{\la}$ in \eqref{e3.15}.

For the proof of the lower estimates for $\wt{A}_{2,m}$, we estimate
from below the factor $r_{\la}$ in \eqref{e3.15}. Since
$-16>-8(2\la+3)$, replacement of $-16$ by $-8(2\la+3)$ in the second
line of \eqref{e3.17} yields
\[
\begin{split}
r_{\la}(m)&\geq (2\la\!+\!3)\big(12m^4\!+\!24\la
m^3\!+\!(12\la^2\!-\!4\la\!-\!48)m^2
\!-\!(4\la^2\!+\!32\la\!-\!8)m\!-\!4\la^2\!-\!4\la\!+\!29\big)\\
&=:(2\la+3)\wt{s}_{\la}(m)\,.
\end{split}
\]

Next, we estimate $\wt{s}_{\la}$ from below, distinguishing between
the cases $-\frac{1}{2}<\la\leq 0$ and  $\la\geq 0$.

If $-\frac{1}{2}<\la\leq 0$, then from
\[
\begin{split}
\wt{s}_{\la}(m)&=12\big(m^2+\la m-\frac{1}{2}\big) \big(m^2+\la
m-\frac{\la}{3}-\frac{7}{2}\big)+8(2\la+1)m-4\la^2-6\la+8\\
&> 12\big(m^2+\la m-\frac{1}{2}\big)\big(m^2+\la
m-\frac{\la}{3}-\frac{7}{2}\big)\,,\quad -\frac{1}{2}<\la\leq 0
\end{split}
\]
(clearly, in that case $8(2\la+1)m-4\la^2-6\la+8>0$) we deduce the
lower bound (i).\smallskip

If $\la\geq 0$, then the lower bound (ii) follows from
\[
\begin{split}
\wt{s}_{\la}(m)&=12\big(m^2+\la m-\frac{1}{2}\big) \big(m^2+\la
m-\frac{\la}{2}-\frac{7}{2}\big)+2\la m^2+(2\la^2+16\la+8)m-4\la^2-7\la+8\\
&> 12\big(m^2+\la m-\frac{1}{2}\big)\big(m^2+\la
m-\frac{\la}{2}-\frac{7}{2}\big)\,,\quad \la\geq 0\,,\; m\geq 2\,.
\end{split}
\]
For the last inequality we have used that, for $\la\geq 0$ and
$m\geq 2$,
$$
g_1(m):=2\la m^2+(2\la^2+16\la+8)m-4\la^2-7\la+8\geq
g_1(2)=33\la+24>0\,.
$$
Lemma \ref{l3.5} is proved. \qed

\section{Estimates for the best Markov constant
\boldmath{$c_n(\la)$}}
Let us recall that $Q_m$ and $\wt{Q}_m$ are the characteristic
polynomials of the matrices $\mathbf{B}_m=\mathbf{A}_m^{-1}$ and
$\wt{\mathbf{B}}_m=\wt{\mathbf{A}}_m^{-1}$, respectively, normalized
by $Q_m(0)=\wt{Q}_m(0)=1$. Hence, their reciprocal polynomials,
\begin{eqnarray}
&&R_m(x)=x^{m}Q_{m}(x^{-1})=x^m-A_{1,m}x^{m-1}+A_{2,m}x^{m-2}-\cdots
+(-1)^m A_{m,m}\,, \label{e4.1}\\
&&\wt{R}_m(x)=x^{m}\wt{Q}_{m}(x^{-1})=x^m-\wt{A}_{1,m}x^{m-1}+\wt{A}_{2,m}x^{m-2}
-\cdots +(-1)^m \wt{A}_{m,m}\,, \label{e4.2}
\end{eqnarray}
are the monic characteristic polynomials of matrices $\mathbf{A}_m$
and $\wt{\mathbf{A}}_m$, respectively. In Sect. 2 we showed that
$Q_m$ and $\wt{Q}_m$ are polynomials orthogonal with respect to
measures supported on the positive axis, therefore their zeros are
single and positive. Then the same observation applies to the zeros
of $R_m$ and $\wt{R}_m$, which we denote by $\{\nu_i\}$ and
$\{\wt{\nu}_i\}$, respectively, so that
\begin{eqnarray*}
&& R_m(x)=(x-\nu_1)(x-\nu_2)\cdots(x-\nu_m)\,,\qquad
0<\nu_1<\nu_2<\cdots<\nu_m\,,\\
&&
\wt{R}_m(x)=(x-\wt{\nu}_1)(x-\wt{\nu}_2)\cdots(x-\wt{\nu}_m)\,,\qquad
0<\wt{\nu}_1<\wt{\nu}_2<\cdots<\wt{\nu}_m\,.
\end{eqnarray*}

Our tool for obtaining two-sided estimates for $\nu_m$ and
$\wt{\nu}_m$ is the following simple observation:
\begin{proposition}\label{p4.1}
Let
$$
f(x)=x^m-a_{1,m}\,x^{m-1}+a_{2,m}\,x^{m-2}-\cdots+(-1)^m\, a_{m,m}
$$
be a polynomial having only real and positive zeros $\{x_i\}$,
$0<x_1\leq x_2\leq\cdots\leq x_m$\,. Then
$$
a_{1,m}-2\,\frac{a_{2,m}}{a_{1,m}}\leq x_m\leq
\sqrt{a_{1,m}^2-2a_{2,m}}\,.
$$
In either place, the equality holds if and only if
$x_1=x_2=\cdots=x_m$.
\end{proposition}

\proof The claim is equivalent to
$$
\frac{x_1^2+x_2^2+\cdots+x_m^2}{x_1+x_2+\cdots+x_m}\leq x_m \leq
(x_1^2+x_2^2+\cdots+x_m^2)^{\frac{1}{2}}\,,
$$
and both the inequalities and the equality cases are obvious. \qed
\smallskip

We obtain separately estimates for $c_n(\la)$ for even and odd $n$.
Theorem~\ref{t1.1} is then obtained as a summary of these results.

\subsection{The cases of even and odd \boldmath{$n$}}

According to Theorem \ref{t2.1}, for the best Markov constant
$c_n(\la)$ we have
\begin{eqnarray}
&&c_{2m}^2(\la)=4\nu_m\,, \label{e4.3}\\
&&c_{2m-1}^2(\la)=4\wt{\nu}_m\,.\label{e4.4}
\end{eqnarray}

\begin{theorem}\label{t4.2}
For all even $n$, $n\geq 4$, and for every $\la>-\Frac{1}{2}$ the
best Markov constant $c_{n}(\la)$ admits the estimates
$$
\frac{(n+2)(n+2\la)(n+\la+\frac{1}{2})^2}{(2\la+1)(2\la+5)} \leq
c_n(\la)^2\leq
\frac{n(n+2\la)(n+2\la+2)\sqrt{(n+2)(n+2\la+3)}}
{2(2\la+1)\sqrt{2\la+5}}
\,.
$$
\end{theorem}

\proof Let us set $n=2m$. We apply Proposition \ref{p4.1} with
$f=R_m$, making use of Lemma~\ref{l3.1}(i) and Lemma~\ref{l3.4}.
\smallskip

1) To derive the lower bound for $c_n(\la)^2$, we estimate
\[
\begin{split}
\nu_m\geq A_{1,m}-2\,\frac{A_{2,m}}{A_{1,m}}&\geq
\frac{m(m+1)(m+\la)(m+\la+1)}{2\la+1}-
\frac{(m-1)(m+1)(m+\la)(m+\la+2)}{2\la+5}\\
&=\frac{(m+1)(m+\la)}{(2\la+1)(2\la+5)}\,
\big[4m(m+\la+1)+(2\la+1)(\la+2)\big]\\
&\geq \frac{(m+1)(m+\la)\big(2m+\la+\frac{1}{2}\big)^2}
{(2\la+1)(2\la+5)}\,.
\end{split}
\]
Hence,
$$
c_{n}^2(\la)=4\nu_m\geq
\frac{4(m+1)(m+\la)\big(2m+\la+\frac{1}{2}\big)^2}
{(2\la+1)(2\la+5)}
=\frac{(n+2)(n+2\la)\big(n+\la+\frac{1}{2}\big)^2}
{(2\la+1)(2\la+5)}\,.
$$

2) For the upper estimate in Theorem~\ref{t4.2}, we have
\[
\begin{split}
\nu_m^2&\leq
A_{1,m}^2-2A_{2,m}\leq\frac{m^2(m+1)^2(m+\la)^2(m+\la+1)^2}{(2\la+1)^2}-
\frac{(m-1)m^2(m+1)(m+\la)^2(m+\la+1)^2}{(2\la+1)(2\la+5)}\\
&=\frac{4m^2(m+\la)^2(m+\la+1)^2(m+1)\big(m+\la+\frac{3}{2}\big)}
{(2\la+1)^2(2\la+5)}=
\frac{n^2(n+2)(n+2\la)^2(n+2\la+2)^2(n+2\la+3)}
{64(2\la+1)^2(2\la+5)}
\end{split}
\]
and then \eqref{e4.3} yields
$$
c_{n}^2(\la)=4\nu_m\leq \frac{n(n+2\la)(n+2\la+2)
\sqrt{(n+2)(n+2\la+3)}} {2(2\la+1)\sqrt{2\la+5}}\,.
$$
The proof of Theorem \ref{t4.2} is complete. \qed

\begin{remark}\label{r4.2}
For $\la\geq 2$ the upper bound for $c_{n}^2(\la)$ in
Theorem~\ref{t4.2} admits a slight improvement, namely, we have
\begin{equation}\label{e4.5}
c_n^2(\la)\leq\frac{n(n+2\la)(n+2\la+2)\sqrt{(n+2)(n+2\la+2)}}
{2(2\la+1)\sqrt{2\la+5}}\,,\qquad \la\geq 2\,.
\end{equation}

Indeed, for $\la\geq 2$ we can replace the lower bound for $A_{2,m}$
in Lemma~\ref{l3.4} by the sharper one
$$
A_{2,m}\geq
\frac{(m-1)m^2(m+1)(m+\la)(m+\la+1)^2(m+\la+2)}{2(2\la+1)(2\la+5)}\,,
$$
and then, proceeding in the same way as above, we arrive at the
estimate \eqref{e4.5}
\end{remark}

\begin{theorem}\label{t4.4}
For all odd $n$, $n\geq 3$, and for every $\la>-\Frac{1}{2}$, the
best Markov constant $c_{n}(\la)$ admits the estimates
$$
\frac{(n+1)\big(n+\la+\frac{1}{2}\big)^2(n+2\la+1)}
{(2\la+1)(2\la+5)}\leq c_n^2(\la)\leq
\frac{(n+1)^{\frac{3}{2}}(n+2\la+1)^{2}(n+2\la'+1)^{\frac{1}{2}}}
{2(2\la+1)\sqrt{2\la+5}}\,,
$$
where $\la'=\max\{\la,0\}$.
\end{theorem}

\proof Let us set $n=2m-1$, $m\geq 2$. We apply
Proposition~\ref{p4.1} with $f=\wt{R}_m$, making use of
Lemma~\ref{l3.1}(ii) and Lemma~\ref{l3.5}.\smallskip

1) For the lower bound, we estimate $\wt{\nu}_m$ from below, using
Proposition~\ref{p4.1}, Lemma~\ref{l3.1}(ii) and Lemma~\ref{l3.5} to
obtain
\[
\begin{split}
\wt{\nu}_m&\geq \wt{A}_{1,m}-2\,\frac{\wt{A}_{2,m}}{\wt{A}_{1,m}}
\geq \frac{m(m+\la)\big(m^2+\la m-\frac{1}{2}\big)}{2\la+1}-
\frac{(m-1)m(m+\la)(m+\la+1)}{2\la+5}\\
&=\frac{m(m+\la)}{(2\la+1)(2\la+5)}\,\Big[4m^2+4\la
m+2\la^2+2\la-\frac{3}{2}\Big]\geq
\frac{m(m+\la)\big(2m+\la-\Frac{1}{2}\big)^2}{(2\la+1)(2\la+5)}\,.
\end{split}
\]
Now the lower estimate for $c_n^2(\la)$ follows from  \eqref{e4.4}:
$$
c_n^2(\la)=4\wt{\nu}_m\geq
\frac{2m(2m+2\la)\big(2m+\la-\Frac{1}{2}\big)^2}{(2\la+1)(2\la+5)}
=\frac{(n+1)(n+2\la+1)\big(n+\la+\Frac{1}{2}\big)^2}{(2\la+1)(2\la+5)}\,.
$$

2)  Next, we prove the upper estimate for $c_n^2(\la)$.\smallskip

2.1) In the case $-\frac{1}{2}<\la\leq 0$, we apply
Proposition~\ref{p4.1}, Lemma~\ref{l3.1}(ii) and inequality (i) in
Lemma~\ref{l3.5} to estimate $\wt{\nu}_m^2$ from above as follows:
{\small
\[
\begin{split}
\wt{\nu}_m^2&\leq \wt{A}_{1,m}^2\!-\!2\wt{A}_{2,m}\leq
\frac{m^2(m\!+\!\la)^2\big(m^2\!+\!\la
m\!-\!\frac{1}{2}\big)^2}{(2\la+1)^2}\!
-\!\frac{(m\!-\!1)m(m\!+\!\la)(m\!+\!\la\!+\!1)\big(m^2\!+\!\la m
\!-\!\frac{1}{2}\big)
\big(m^2\!+\!\la m\!-\!\frac{\la}{3}\!-\!
\frac{7}{2}\big)}{(2\la+1)(2\la+5)}\\
&=\frac{m(m\!+\!\la)\big(m^2\!+\!\la m
\!-\!\frac{1}{2}\big)}{(2\la+1)^2(2\la+5)}\,\Big[
(2\la\!+\!5)m(m\!+\!\la)\big(m^2\!+\!\la m\!-\!\frac{1}{2}\big)
\!-\!(2\la\!+\!1)(m\!-\!1)(m\!+\!\la\!+\!1)\big(m^2\!+\!\la
m\!-\!\frac{\la}{3}\!-\!\frac{7}{2}\big)\Big]\\
&=\frac{4m(m\!+\!\la)\big(m^2\!+\!\la m
\!-\!\frac{1}{2}\big)}{(2\la+1)^2(2\la+5)}\,\Big[m^2(m+\la)^2 \!+\!
\Big(\frac{2}{3}\la^2\!+\!\frac{7}{3}\la\!+\!\frac{1}{2}\Big)m(m\!+\!\la)
\!-\!\frac{1}{4}(2\la\!+\!1)(\la\!+\!1)\Big(\frac{\la}{3}\!+\!\frac{7}{2}\Big)\Big]\\
&\leq \frac{4m^2(m\!+\!\la)^2\big(m^2\!+\!\la m
\!-\!\frac{1}{2}\big)}{(2\la+1)^2(2\la+5)}\,\Big[m(m+\la)\!+\!
\frac{2}{3}\la^2\!+\!\frac{7}{3}\la\!+\!\frac{1}{2}\Big]\,.
\end{split}
\]}

Since $g_2(\la):=\Frac{2}{3}\la^2+\Frac{7}{3}\la+\frac{1}{2}$ is a
monotone increasing function in $(-1/2,0]$, the expression in the
last brackets does not exceed $m^2+\la m+\Frac{1}{2}$, hence
$$
\wt{\nu}_m^2\leq\frac{4m^2(m+\la)^2\big[m^2(m+\la)^2-\frac{1}{4}\big]}
{(2\la+1)^2(2\la+5)}\leq \frac{4m^4(m+\la)^4}{(2\la+1)^2(2\la+5)}\,.
$$
Now from \eqref{e4.4} we obtain the desired upper estimate for
$c_n^2(\la)$:
$$
c_n^2(\la)=4\wt{\nu}_m\leq\frac{8m^2(m+\la)^2}{(2\la+1)^2\sqrt{2\la+5}}
=\frac{(n+1)^2(n+2\la+1)^2}{2(2\la+1)\sqrt{2\la+5}}=
\frac{(n+1)^{\frac{3}{2}}(n+2\la'+1)^{\frac{1}{2}}(n+2\la+1)^{2}}
{(2\la+1)\sqrt{2\la+5}}\,.
$$

2.2) In view of \eqref{e4.4}, in the case $\la\geq 0$ the upper
estimate for $c_n^2(\la)$ in Theorem~\ref{t4.4} is equivalent to
\begin{equation}\label{e4.6}
\wt{\nu}_m^2\leq \frac{4m^3(m+\la)^5}{(2\la+1)^2(2\la+5)}\,.
\end{equation}

We apply Proposition~\ref{p4.1}, Lemma~\ref{l3.1}(ii) and inequality
(ii) in Lemma~\ref{l3.5} to estimate $\wt{\nu}_m^2$ from above as
follows: {\small
\[
\begin{split}
\wt{\nu}_m^2\leq&\frac{m^2(m+\la)^2\big(m^2+\la
m-\frac{1}{2}\big)^2}{(2\la+1)^2}
-\frac{(m-1)m(m+\la)(m+\la+1)\big(m^2+\la
m-\frac{1}{2}\big)\big(m^2+\la
m-\frac{\la+7}{2}\big)}{(2\la+1)(2\la+5)}\\
=&\frac{m(m\!+\!\la)\big(m^2\!+\!\la
m\!-\!\frac{1}{2}\big)}{(2\la+1)^2(2\la+5)}\,
\Big[(2\la\!+\!5)(m^2\!+\la m) \big(m^2\!+\!\la
m\!-\!\frac{1}{2}\big) \!-\!(2\la\!+\!1)(m^2\!+\!\la m
\!-\!\la\!-\!1)\big(m^2\!+\!\la m\!-\!\frac{\la\!+\!7}{2}\big)\Big]\\
=&\frac{m(m\!+\!\la)\big(m^2\!+\!\la
m\!-\!\frac{1}{2}\big)}{(2\la+1)^2(2\la+5)}\,\Big[4m^2(m\!+\!\la)^2
\!+\!\frac{1}{2}(6\la^2+19\la+4)m(m\!+\!\la)
\!-\!\frac{1}{2}(2\la\!+\!1)(\la\!+\!1)(\la\!+\!7)\Big]\\
\leq&\frac{4m^2(m\!+\!\la)^2}{(2\la+1)^2(2\la+5)}\,\big(m^2\!+\!\la
m\!-\!\frac{1}{2}\big)\Big[m^2\!+\la
m+\frac{1}{8}(6\la^2+19\la+4)\Big]\,.
\end{split}
\]}

To prove \eqref{e4.6}, it suffices to show that
$$
\big(m^2+\la m-\frac{1}{2}\big)\Big[m^2+\la
m+\frac{1}{8}(6\la^2+19\la+4)\Big]\leq m(m+\la)^3\,, \qquad \la\geq
0\,,\ \ m\geq 2\,.
$$

For $m\geq 3$ the above inequality follows from
$$
m(m+\la)^3-\big(m^2+\la m-\frac{1}{2}\big)\Big[m^2+\la
m+\frac{1}{8}(6\la^2+19\la+4)\Big]=\frac{1}{8}\,\la m
(m+\la)(8m+2\la-19)+\frac{1}{16}(6\la^2+19\la+4)\,,
$$
while for $m=2$ it is equivalent to the inequality
$$
8\la^3+10\la^2-5\la+4\geq 0\,,\quad \la\geq 0,
$$
which is readily verified to be true. \qed

\subsection{Proof of Theorem 1.1 and Corollary 1.3}

\textbf{Proof of Theorem 1.1.} Clearly, the lower bound for
$c_n^2(\la)$ in Theorem~\ref{t1.1} is smaller than the lower bounds
in Theorems~\ref{t4.2} and \ref{t4.4}, hence it is a lower bound in
the cases of both even and odd $n$. \smallskip

Next, we prove the upper bound for $c_n^2(\la)$ in
Theorem~\ref{t1.1}. To this end, we apply the geometric mean -
arithmetic mean inequality in Theorems~\ref{t4.2} and \ref{t4.4}  to
obtain
$$
c_n^2(\la)\leq \Frac{\big(
n+\Frac{5}{4}\la+\Frac{9}{8}\big)^4}{2(2\la+1)\sqrt{2\la+5}}\,,\qquad
n=2m\,,
$$
$$
c_n^2(\la)\leq \Frac{\big(
n+\la+\Frac{\la'}{4}+1\big)^4}{2(2\la+1)\sqrt{2\la+5}}\,,\qquad
n=2m-1\,,\ \la'=\max\{0,\la\}
$$
and compare the right-hand sides of these inequalities, observing
that the first one is the greater. \qed

\begin{remark}\label{r4.5}
Applying the geometric mean -- arithmetic mean inequality to the
upper bounds for $c_n^2(\la)$ in Theorems~\ref{t4.2} and \ref{t4.4}
to obtain the upper bound in Theorem~\ref{t1.1}, we certainly lose.
For instance, for a fixed $n$, the upper bounds in
Theorems~\ref{t4.2} and \ref{t4.4} are $O(\la)$ as
$\la\rightarrow\infty$ (notice that the same applies to the lower
bounds therein!), while the resulting upper bound in
Theorem~\ref{t1.1} is $O(\la^{\frac{5}{2}})$ as
$\la\rightarrow\infty$. However, as was already said, the upper
estimates here are good for relatively small $\la$, say, $\la\leq
25$. For big $\la$, we have the better upper estimates \eqref{e2} in
Theorem~A.
\end{remark}

\noindent \textbf{Proof of Corollary 1.3.} The comparison of the
bounds for $c_n^2(\la)$ in Theorems~\ref{t4.2} and \ref{t4.4}
reveals that for $\la<\Frac{1}{2}$ the smaller lower bound is the
one in Theorem~\ref{t4.2}, while in the limit case $\la=-1/2$ the
bigger numerator has the upper bound in Theorem~\ref{t4.4}. By
taking the limits in the expressions obtained from corresponding
bounds we obtain the result. \qed\bigskip

\noindent \textbf{Acknowledgements} This research was supported by
the Bulgarian National Science Fund under Contract DN 02/14.

\bigskip

\noindent
{\sc Dragomir Aleksov, Geno Nikolov} \smallskip\\
Department of Mathematics and Informatics\\
Universlty of Sofia \\
5 James Bourchier Blvd. \\
1164 Sofia \\
BULGARIA \\
{\it E-mails:} {\tt dragomira@fmi.uni-sofia.bg,
geno@fmi.uni-sofia.bg}

\end{document}